\documentclass[11pt]{article}
\usepackage{amsmath, amssymb, physics}
\usepackage{amscd}
\usepackage{geometry}
\usepackage{graphicx}
\usepackage{subfigure}
\usepackage{xcolor}
\usepackage{hyperref}
\usepackage{forest}
\usepackage{epsfig}
\usepackage{tikz}
\usetikzlibrary{positioning}

\forestset{
  */.style={
    delay+={append={[]},}
  },
  rooted tree/.style={
    for tree={
      grow'=90,
      parent anchor=center,
      child anchor=center,
      s sep=2.5pt,
      if level=0{
        baseline
      }{},
      delay={
        if content={*}{
          content=,
          append={[]}
        }{}
      }
    },
    before typesetting nodes={
      for tree={
        circle,
        fill,
        minimum width=3pt,
        inner sep=0pt,
        child anchor=center,
      },
    },
    before computing xy={
      for tree={
        l=5pt,
      }
    }
  },
  big rooted tree/.style={
    for tree={
      grow'=90,
      parent anchor=center,
      child anchor=center,
      s sep=2.5pt,
      if level=0{
        baseline
      }{},
      delay={
        if content={*}{
          content=,
          append={[]}
        }{}
      }
    },
    edge={
        semithick,
        -Latex
    },
    before typesetting nodes={
      for tree={
        circle,
        fill,
        minimum width=5pt,
        inner sep=0pt,
        child anchor=center,
      },
    },
    before computing xy={
      for tree={
        l=15pt,
      }
    }
  },
  my node/.style={
    circle, fill=#1,
  },
  gray node/.style={
    my node=gray, inner sep=#1mm,
  },
  red node/.style={
    my node=red, inner sep=#1mm,
  },
  blue node/.style={
    my node=blue, inner sep=#1mm,
  },
  green node/.style={
    my node=green, inner sep=#1mm,
  },
}

\DeclareDocumentCommand\rt{o}{\Forest{rooted tree [#1]}}
\DeclareDocumentCommand\bigrt{o}{\Forest{big rooted tree [#1]}}

\usepackage{amsthm}

\newtheorem{thm}{Theorem}

\newtheorem{corollary}{Corollary}
\newtheorem{prop}{Proposition}

\let\emptyset\varnothing

\title{Hopf Algebra of Graphs for Operator Splitting Methods}
\author{Gabriel Barrios de Leon}

\begin{document}

\maketitle

\begin{abstract}
    Motivated by combinatorial identities arising in the study of Operator Splitting methods and their description in terms of bipartite rooted trees, we establish a connection between Operator Splitting methods and the theory of combinatorial Hopf algebras. We show that the inclusion--exclusion expansion of the elementary weight of a rooted tree over its rooted subforests can be interpreted as the antipode of Schmitt's Hopf algebra of graphs. This interpretation recovers known identities for operator splitting and partitioned Runge–Kutta methods from general properties of the antipode. We further show that adjoint operator splitting methods can be described algebraically through the antipode and the canonical involution of the original method. The resulting framework also yields new combinatorial identities and connects the analysis of geometric numerical integrators with the theory of incidence Hopf algebras.
\end{abstract}

\section{Introduction}
Graph theory has played a fundamental role in the development of modern numerical time integration methods by providing a combinatorial framework for analyzing the order conditions of high-order integrators. A fundamental contribution in this direction was the pioneering work of Butcher, who observed that the elementary differentials and elementary weights associated with Runge-Kutta (RK) methods can be naturally indexed by rooted trees \cite{butcher1963coefficients,butcher1972algebraic}. This correspondence transformed the derivation of order conditions into a combinatorial problem, allowing the structure of rooted trees to be exploited in the systematic construction and analysis of numerical integrators. Since then, rooted trees have become an indispensable tool in geometric numerical integration, giving rise to the theory of B-series and revealing deep algebraic structures underlying RK and related methods \cite{hairer2006geometric,hairer2003geometric,butcher2021bseries}.

The combinatorial structure of B-series through rooted trees led to the discovery of an algebraic structure underlying numerical integration methods. In particular, the set of rooted trees can be endowed with the structure of a Hopf algebra \cite{chartier2011hopf,brouder2000runge}, providing a powerful framework for studying the composition, substitution, and inversion of B-series. Remarkably, this Hopf algebra was independently discovered in the context of perturbative quantum field theory by Connes and Kreimer \cite{connes1998hopf}, in the context of renormalization. 

Operator splitting (OS) methods, on the other hand, constitute an important class of numerical integrators for initial value problems whose vector field can be decomposed into a sum of simpler components \cite{hairer2003geometric}. These methods construct an approximation by composing the exact or approximate flows associated with the individual subproblems.   As with RK methods, the design and analysis of high-order splitting schemes rely on understanding intricate algebraic and combinatorial order conditions, motivating the search for mathematical structures that organize these conditions in a systematic way. Traditionally, these order conditions are derived through repeated applications of the Baker--Campbell--Hausdorff (BCH) formula \cite{blanes2016introduction,auzinger2016bch}, which expresses the composition of exponentials of non-commuting operators as the exponential of an infinite series of nested commutators. By recursively applying the BCH formula to the composition defining a splitting method and comparing the resulting Lie series with the exact flow, one obtains the algebraic constraints on the coefficients required to achieve a prescribed order of accuracy.

The present work is motivated by the study of fractional step Runge--Kutta (FSRK) methods through the combinatorial properties of the elementary weights associated to bipartite rooted trees \cite{barrios2026FSRK,barrios2026splitting}. In that setting, elementary weights can be expressed in terms of the operator-splitting coefficients and expanded over subforests obtained by deleting edges from the original tree. These expansions lead to identities involving the tree factorial and reveal a systematic pattern of alternating signs. The recurrence of such identities naturally suggested the existence of an underlying algebraic framework governing these combinatorial phenomena. 

In this paper we develop this connection using Schmitt's Hopf algebra of graphs for the case of trees, where the expansions over subforests and the associated inversion formulas become manifestations of the coproduct and antipode. This viewpoint not only explains and unifies the identities obtained for FSRK methods but also extends them to a broader class of OS schemes. We first formulate the subforest expansions and the associated identities for elementary weights, including their interpretation through the edge-deletion poset and Möbius inversion. We then show how the expansion operator is related to the antipode, providing an algebraic explanation for identities that arise in the analysis of FSRK and operator-splitting methods. Finally, we apply this framework to operator splitting—including adjoint and palindromic methods—and to partitioned Runge–Kutta methods. In addition to recovering known identities, these applications yield new combinatorial relations and clarify the role of incidence-algebra structures in geometric numerical integration.

The paper is organized as follows. In section \ref{sec:background} we recall the basics of graph theory and introduce in Theorems \ref{thm:phiC} and \ref{thm:phi_homeormorphic_candycanes} the two identities relating the characters $\phi$ and $\phi^C$, which are the elementary weights for bipartite trees in the context of 2-operator splitting FSRK methods. We recognize that labeled rooted trees to study the BCH Formula have already been used (see, for example, \cite{eizinger2009antipode,iacono2013posetted}). Section \ref{sec:incidence} introduces the edge-deletion poset, and we show that the characters $\phi$ and $\phi^C$ are related by Möbius inversion. 
Section \ref{sec:hopf_algebra} develops the necessary Hopf-algebraic framework and identifies the subforest expansion with the antipode of Schmitt’s Hopf algebra of graphs \cite{SCHMITT1994incidenceHopf,humpert2012incidence,denneberg2012incidence}. Finally, Section \ref{sec:applications} is devoted to illustrate applications of the theory developed in this work for numerical integration, in particular, for OS methods (including adjoint and palindrome methods) and partitioned Runge-Kutta methods.

\section{Bipartite trees in FSRK methods}
\label{sec:background}
Consider the initial value problem
\begin{equation}
u'(t)=f(u,t), \qquad u(0)=u_0,
\label{eq:ivp}
\end{equation}
and assume that the vector field can be decomposed as
\begin{equation}
f(u,t)=\sum_{\ell=1}^{N}f^{\ell}(u,t).
\label{eq:ivp_fsrk}
\end{equation}
Such decompositions naturally arise in many applications (for example, reaction-diffusion systems, fluid dynamics, Hamiltonian mechanics, etc.), where each operator $f^{\ell}$ models a different physical phenomenon or corresponds to a distinct term in the governing equations. Rather than treating the complete system as a whole, one may exploit this decomposition to construct numerical methods that separately approximate the evolution generated by each component.

OS methods are based precisely on this principle. They approximate the solution of \eqref{eq:ivp} by composing the solution operators associated with the individual subproblems, thereby replacing the original evolution by a sequence of simpler ones. This strategy often leads to computationally efficient schemes and, in many situations, inherits desirable geometric or physical properties from the constituent flows. Among the most widely used examples are the Lie--Trotter and Strang splittings, which provide first and second-order approximations, respectively \cite{blanes2024splitting}.

FSRK methods build upon this idea by replacing the exact evolution of each subproblem with a suitable RK discretization \cite{hairer2003geometric}. Then, the overall method combines the flexibility of RK schemes with the structure of OS methods, making it possible to construct high-order integrators while preserving the advantages associated with the decomposition of the original problem. 

\subsection{Preliminaries on rooted trees}
Recall that a graph $G$ is a pair $(V(G),E(G))$, where $V(G)$ and $E(G)$ denote the sets of vertices and edges, respectively. The number of edges coming out from a vertex is called the \textit{degree} of the vertex. A graph containing no cycles is called a \emph{forest}, while a connected forest is called a \emph{tree}. Of particular importance in numerical analysis are \emph{rooted trees}, namely trees in which one distinguished vertex, called the \emph{root}, is specified. The choice of the root induces a natural orientation on the edges, directed away from the root, and thereby establishes a hierarchical structure on the vertices.

For every non-root vertex $v$, there exists a unique adjacent vertex lying on the path from $v$ to the root. This vertex is called the \emph{parent} of $v$, and $v$ is said to be a \emph{child} of its parent. Vertices having no children are referred to as \emph{leaves}. Finally, the order of a rooted tree $t$, denoted by $|t|$, is defined as the number of its vertices. Also, we use the notation $t = [t_1,\dots,t_m]$ to denote a tree with subtrees $t_i,\; i=1,\dots,m$. In case we are dealing with colored trees, we use the notation $t = [t_1,\dots,t_m]_{\ell}$ to indicate that $t$ has root with color $\ell$.

\subsection{ARK representation}
Using the additive Runge-Kutta (ARK) representation of FSRK methods developed in \cite{Spitieri2023FSRK}, we studied the construction of higher order methods with $N=2$ in \cite{barrios2026FSRK}. Subsequently, in \cite{barrios2026splitting}, we studied  FSRK methods with $N=2$ constructed from $s-$stages OS methods with coefficients $\{\alpha_i^{\ell} \}$, $i = 1,\dots,s$, $\ell = \mu,\nu$. There we focused on the elementary weights associated with \emph{bipartite rooted trees}, namely rooted trees whose vertices are assigned one of two colors in such a way that the endpoints of every edge have different colors. Among the characters defined on rooted trees, the \textbf{tree factorial} (or density) is particularly useful for numerical analysis:
\begin{equation}
        t! = [t_1,\dots,t_m]!= 
        \begin{cases}
            1, & |t| =1,\\
            |t|\cdot t_1! \cdots t_m!, & \text{else}.
        \end{cases}
\end{equation}
Finally, we define an equivalence relation for two rooted trees $t_1 ,t_2$ so that $t_1 \sim t_2$ if they have the same structure as unrooted trees, and differ only on the choice of the vertex that is the root.

Let $T_{cc}^p$ be the set of bipartite trees of order $p$, and let $T_{cc}$ be the set of all bipartite rooted trees; let $\{\alpha_i^{\ell} \}$ be the OS coefficients of an FSRK, given in its ARK representation, built with RK subintegrators of (at least) first order of accuracy. We showed that the elementary weights associated to bipartite trees, $\phi(t):T_{cc} \rightarrow \mathbb{R}$, are given recursively by: 
    \begin{equation}
        \phi(t) =
        \begin{cases}
            1, &\quad |t| = 1,\\
            \sum_{i=1}^s\alpha(t)_i\prod_{j=1}^m\sum_{k=1}^{i-\delta(t_j)}\tilde \phi_k(t_j), &\quad t=[t_1,\dots,t_m]_\ell,
        \end{cases}
        \label{eq:phi_def}
    \end{equation}
    where
    \begin{equation}
    \tilde \phi_i(t) =
    \begin{cases}
        \alpha(t)_i, &\quad |t| = 1,\\
        \alpha(t)_i \prod_{j=1}^m\sum_{k=1}^{i-\delta(t_j)}\tilde \phi_k(t_j), &\quad t= [t_1,\dots,t_m]_\ell,
    \end{cases}
    \label{eq:phi_tilde}
    \end{equation}
    and 
    \begin{equation}
        \delta(t=[t_1,\dots,t_m]_\ell) = 
        \begin{cases}
            1, \quad \ell(t) =\nu,\\
            0, \quad \ell(t) = \mu.
        \end{cases}    
    \end{equation}

We assume that the OS methods are (at least) first order, which translates into the coefficients satisfying $\sum_i\alpha_i^{\ell} =1, \; \ell = \mu,\nu $. Then, we defined a \textit{complementary} character $\phi^C(t):T_{cc} \rightarrow \mathbb{R}$ that inverts the order of the OS coefficients:
\begin{equation}
        \phi^C(t) =
        \begin{cases}
            1, &\quad |t| = 1,\\
            \sum_{i=1}^s\alpha(t)_i\prod_{j=1}^m\sum_{k=i + 1- \delta(t_j)}^{s}\tilde \phi^C_k(t_j), &\quad t=[t_1,\dots,t_m]_\ell,
        \end{cases}
        \label{eq:phi_complementary}
    \end{equation}
    where
    \begin{equation}
    \tilde \phi^C_i(t) =
    \begin{cases}
        \alpha(t)_i, &\quad |t| = 1,\\
        \alpha(t)_i \prod_{j=1}^m\sum_{k=i +1- \delta(t_j)}^{s}\tilde \phi^C_k(t_j), &\quad t= [t_1,\dots,t_m]_\ell.
    \end{cases}
    \end{equation}

Using the operation of deleting edges, we found an equation that relates the two characters for a given $t \in T_{cc}$ making use of a linear expansion of all proper subforests of $t$ (i.e., the induced subforest obtained after deleting a subset of edges and keepig the complete set of vertices). For an edge $e \in E(t)$, define $c_e(t)$ as a cut that returns the rooted subforest obtained after deleting the edge $e$ from $t$.  Now, we define the following operator:
\begin{equation}
    K(t) := \prod_{e\in E(t)} (c_e -1)(t).
    \label{eq:K_def}
\end{equation}
For the algebraic interpretation of $K$, it is useful, as will be seen later, to rewrite the cuts in terms of subsets of edges. We extend the operation of deleting single edges to subsets of edges $A \subseteq E(t)$ by $c_A(t)$.
\begin{prop}
    $K$ can be rewritten as
\begin{equation}
    K(t) = \sum_{A\subseteq E(t)}(-1)^{|E(t)|-|A|}c_A(t).
    \label{eq:K_expansion}
\end{equation}
\end{prop}
\noindent
Then, the characters $\phi$ and $\phi^C$ are related by the following identity:
\begin{equation}
    \phi(t) = \sum_{A \subseteq E(t)}(-1)^{|E(t)| - |A|}\phi^C\left( c_A(t)\right)  = \phi^C\left( K(t)\right).
    \label{eq:phi_subgraphs}
\end{equation}
As illustrations, consider $t = \rt[[,color=white,draw=black [][]]]$, for which \eqref{eq:phi_subgraphs} is
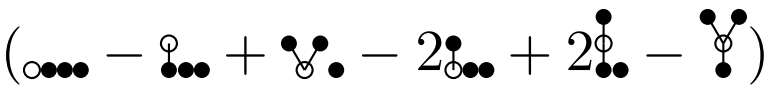
\begin{equation}
    \phi(t) = \phi^C(\rt[,color=white,draw=black]\rt[]\rt[]\rt[] -\rt[[,color=white,draw=black]]\rt[]\rt[] +\rt[,color=white,draw=black[][]]\rt[] - 2\rt[,color=white,draw=black[]]\rt[]\rt[] + 2\rt[[,color=white,draw=black []]]\rt[] - \rt[[,color=white,draw=black [][]]]). 
    \label{eq:example_subforest1}
\end{equation}
For a higher order example, consider the sixth order tree $t=\rt[[,color=white,draw=black[]][,color=white,draw=black[][]]]$.


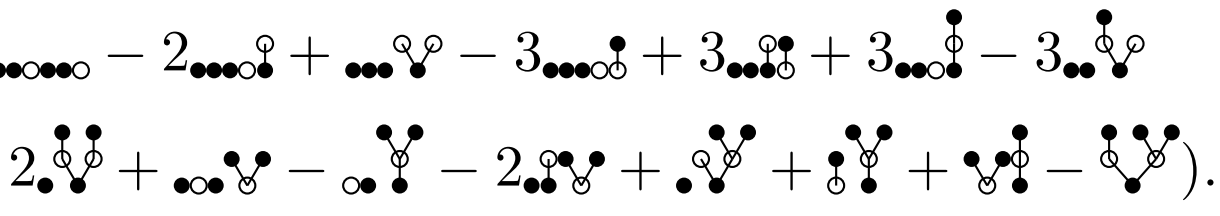
\begin{equation}
    \phi(t) = \phi^C(\rt[]\rt[]\rt[,color=white,draw=black]\rt[]\rt[]\rt[,color=white,draw=black]     - 2\rt[]\rt[]\rt[]\rt[,color=white,draw=black]  \rt[[,color=white,draw=black]]  + \rt[]\rt[]\rt[]  \rt[[,color=white,draw=black][,color=white,draw=black]] - 3\rt[]\rt[]\rt[]\rt[,color=white,draw=black]  \rt[,color=white,draw=black[]] + 3\rt[]\rt[]   \rt[[,color=white,draw=black]]\rt[,color=white,draw=black[]]
    +3\rt[]\rt[]  \rt[,color=white,draw=black] \rt[[,color=white,draw=black[]]] - 3\rt[]\rt[]   \rt[[,color=white,draw=black[]][,color=white,draw=black]] 
\end{equation}
$$
+2\rt[]\rt[]   \rt[,color=white,draw=black[]]\rt[,color=white,draw=black[]] - 4\rt[]  \rt[,color=white,draw=black[]]\rt[[,color=white,draw=black[]]] +2\rt[]  \rt[[,color=white,draw=black[]][,color=white,draw=black[]]]+\rt[]\rt[,color=white,draw=black]\rt[]  \rt[,color=white,draw=black[][]] -\rt[,color=white,draw=black]\rt[]  \rt[[,color=white,draw=black[][]]] - 2\rt[]  \rt[[,color=white,draw=black]]\rt[,color=white,draw=black[][]] + \rt[]  \rt[[,color=white,draw=black][,color=white,draw=black[][]]] +     \rt[,color=white,draw=black[]]\rt[[,color=white,draw=black[][]]] +     \rt[,color=white,draw=black[][]]\rt[[,color=white,draw=black[]]] -     \rt[[,color=white,draw=black[]][,color=white,draw=black[][]]]).
$$

We found a sufficient (but not necessary) condition for the two character to be equivalent. 
\begin{thm}[Theorem 4, \cite{barrios2026splitting}]\label{thm:phiC}
    Let $\phi(t) := \frac{1}{t!}$ for all bipartite trees such that $|t| <p$, and assume that for some $t_p \in T^p_{cc}$ we set $\phi(t_p) = \frac{1}{t_p!}$. Then, $\phi^C(t_p) = \frac{1}{t_p!} =\phi(t_p)$.
\end{thm}
As an example, consider again the tree in \eqref{eq:example_subforest1}, which yields
$$
\frac{1}{12} = 1 - \frac{1}{2} + \frac{1}{3} -1 + \frac{2}{6} - \frac{1}{12}.
$$
In fact, the elementary weights associated to bipartite trees belonging to the same equivalence class $\sim$ are equivalent.
\begin{thm}[Theorem 5, \cite{barrios2026splitting}]
    Assume $\phi(t) := \frac{1}{t!}$ for all bipartite rooted trees such that $|t| <p$, and let $t_p,t_p' \in T_{cc}^p$  be such that $t_p\sim t_p'$. If for some $t_p \in T^p_{cc}$ we set $\phi(t_p) = \frac{1}{t_p!}$, then $\phi(t_p') = \frac{1}{t_p'!} $.
    \label{thm:phi_homeormorphic_candycanes}
\end{thm}
\noindent
As a Corollary of Theorem \ref{thm:phi_homeormorphic_candycanes}, we find a combinatorial property of the tree factorial for rooted trees belonging the same equivalence class as unrooted trees.
\begin{corollary}\label{corollary:edge_rooted}
    Let $t$ be an edge-rooted tree. Let the root be the edge $e=(v_1,v_2)$, and let $t_1,t_2$ be the two rooted trees obtained after deleting $e$, rooted in $v_1,v_2$ respectively. Then,
    \begin{equation}
        \frac{1}{(t_1 \triangleleft t_2)!} + \frac{1}{(t_2 \triangleleft t_1)!} = \frac{1}{t_1!}\frac{1}{t_2!}.
        \label{eq:partitioned}
    \end{equation}
\end{corollary}
This identity \eqref{eq:partitioned} is also found in \cite{Abia1993partitioned} in the context of PRK for separable Hamiltonian systems. We can take advantage of this result and extend it as follows.
\begin{prop}\label{prop:degree}
    Let $t$ be a tree with sets $(V(t),E(t))$ of vertices and edges. For any $v\in V(t)$, let $T_v$ denote the rooted tree rooted at $v$, and for any edge $e=(v_1,v_2) \in E(t)$, let $C_e(t)$ be the rooted forest consisting of the two trees rooted at $v_1$ and $v_2$ obtained after deleting $e$. Then, the character $\phi(t) = \frac{1}{t!}$ satisfies
    \begin{equation}
        \sum_{e\in E(t)} \phi \left( C_e(t) \right) = \sum_{v\in V(t)} deg(v)\phi(T_v).
        \label{eq:degree}
    \end{equation}
\end{prop}
\begin{proof}
    Setting $deg(v)$ as the number of edges incident to $v$, using Corollary \ref{corollary:edge_rooted} for each vertex $v \in V(t)$ provides Equation \eqref{eq:degree}.
\end{proof}
\noindent
Before diving into  algebraic structures, we provide a combinatorial property of the operator $K$:
\begin{prop}
\label{prop:involution}
$K$ is an involution:
    \begin{equation}
        K^2 := K\circ K  = id.
    \end{equation}
\end{prop}

\begin{proof}
To prove that $K^2 = \text{id}$, we analyze the action of the operator $K$ on a tree $t$ with an edge set $E(t)$ containing $n$ edges, i.e., $|E(t)| = n$. By definition, the operator expands as an alternating sum over all possible subsets of edges to be deleted:
\[ K(t) = \sum_{A \subseteq E(t)} (-1)^{n - |A|} (t \setminus A) \]
When we apply the operator $K$ a second time, it acts on the remaining forest $t \setminus A$. At this stage, the available edge set has shrunk to $E(t) \setminus A$, which contains $n - |A|$ edges. Applying the definition of $K$ to this new structure yields:
\[ K(t \setminus A) = \sum_{B \subseteq E(t) \setminus A} (-1)^{(n - |A|) - |B|} (t \setminus A \setminus B) \]
Substituting this expression back into the equation for $K^2(t)$ gives a nested double sum over pairs of disjoint edge subsets $A$ and $B$:
\[ K^2(t) = \sum_{A \subseteq E(t)} \sum_{B \subseteq E(t) \setminus A} (-1)^{n - |A|} (-1)^{n - |A| - |B|} (t \setminus (A \cup B)) \]
We can simplify the sign exponent by combining the powers of $-1$:$ (n - |A|) + (n - |A| - |B|) = 2n - 2|A| - |B| $. Since $2n$ and $-2|A|$ are always even integers for any choice of $n$ and $A$, they do not affect the sign. Thus, $(-1)^{2n - 2|A| - |B|} = (-1)^{-|B|} = (-1)^{|B|}$. The expression for $K^2(t)$ simplifies to:
\[ K^2(t) = \sum_{A \subseteq E(t)} \sum_{B \subseteq E(t) \setminus A} (-1)^{|B|} (t \setminus (A \cup B)) \]
To evaluate this sum, we change variables by letting $C = A \cup B$ represent the total set of edges deleted after both operations. Since $A$ and $B$ are disjoint, specifying $C$ and $B$ completely determines $A$ as $C \setminus B$. We can restructure the calculation to choose the final subset $C$ first, and then sum over all possible subsets $B$ contained within $C$:
\[ K^2(t) = \sum_{C \subseteq E(t)} \left( \sum_{B \subseteq C} (-1)^{|B|} \right) (t \setminus C) \]
The inner sum $\sum_{B \subseteq C} (-1)^{|B|}$ is the standard alternating sum over all subsets of $C$. By the binomial theorem, this sum is equal to $(1 - 1)^{|C|}$. This yields two distinct cases depending on the size of $C$:
\[ \sum_{B \subseteq C} (-1)^{|B|} = \begin{cases} 1 & \text{if } C = \emptyset \\ 0 & \text{if } C \neq \emptyset \end{cases} \]
Every term where $C \neq \emptyset$ vanishes due to this cancellation. The only term that survives the entire double sum is the one where $C = \emptyset$, which forces both $A = \emptyset$ and $B = \emptyset$. Evaluating this surviving term gives:
\[ K^2(t) = 1 \cdot (t \setminus \emptyset) = t \]
Therefore, $K^2 = \text{id}$.

\end{proof}

\section{Incidence algebras of posets}
\label{sec:incidence}
% The edge-deletion poset of a rooted tree.
\subsection{The edge-deletion poset of a rooted tree}
\label{sec:edge_deletion}
Theorems \ref{thm:adjoint_character} and \ref{thm:phi_homeormorphic_candycanes} allow us to express the character associated with a rooted tree as a sum over all subforests obtained by deleting edges. Such expansions suggest that the collection of subforests carries a natural partially ordered structure. In this section, we formalize this observation by introducing the edge-deletion poset associated with a rooted tree.

Let $t=(V,E)$ be a rooted tree (not necessarily bipartite nor colored). Since deleting edges cannot create cycles, every such graph is again a forest whose connected components inherit the roots induced by $t$.

We denote by
\[
\mathcal{F}(t)
    =
    \{F \subseteq t:\, F \text{ is obtained by deleting edges from } t\}
\]
the collection of all proper subforests of $t$.

A natural partial order on $\mathcal{F}(t)$ is given by edge inclusion. Namely, for
$F_1,F_2\in\mathcal{F}(t)$, we write
\[
F_1\leq F_2
\]
if every edge of $F_1$ is also an edge of $F_2$, i.e.,
\[
E(F_1)\subseteq E(F_2).
\]
Equivalently, $F_1$ can be obtained from $F_2$ by deleting edges.

It is immediate that $(\mathcal{F}(t),\leq)$ is a finite partially ordered set. The minimal element is the edgeless forest consisting of the isolated vertices of $t$, while the maximal element is the tree $t$ itself.

More generally, given two subforests $F_1\leq F_2$, the interval
\[
[F_1,F_2]
    =
    \{F\in\mathcal{F}(t):F_1\leq F\leq F_2\}
\]
consists of all forests obtained by deleting from $F_2$ an arbitrary subset of the edges belonging to
$E(F_2)\setminus E(F_1)$.
Consequently,
\[
|[F_1,F_2]|
    =
    2^{|E(F_2)|-|E(F_1)|},
\]
so every interval is finite. Hence, $(\mathcal{F}(t),\leq)$ is a locally finite poset. Figure \ref{fig:edge_poset} illustrates the poset of edge-deletion for a fourth order tree.


\begin{figure}[ht]
\centering

%-------------------------------------------------
% Panel (a): Edge deletion poset
%-------------------------------------------------
\subfigure[]{
\begin{tikzpicture}[
    node distance=1.25cm and 2cm,
    every node/.style={
        draw,
        rounded corners,
        inner sep=2pt,
        font=\small
    },
    >=latex
]

\node (top) {$\varnothing$};

\node (e1) [below left=of top] {$\{e_1\}$};
\node (e2) [below=of top] {$\{e_2\}$};
\node (e3) [below right=of top] {$\{e_3\}$};

\node (e12) [below left=of e2] {$\{e_1,e_2\}$};
\node (e13) [below right=of e2] {$\{e_1,e_3\}$};
\node (e23) [below=of e2] {$\{e_2,e_3\}$};

\node (all) [below=of e23] {$\{e_1,e_2,e_3\}$};

\draw (top)--(e1);
\draw (top)--(e2);
\draw (top)--(e3);

\draw (e1)--(e12);
\draw (e1)--(e13);

\draw (e2)--(e12);
\draw (e2)--(e23);

\draw (e3)--(e13);
\draw (e3)--(e23);

\draw (e12)--(all);
\draw (e13)--(all);
\draw (e23)--(all);

\end{tikzpicture}
\label{subfig:subforest_deletion}
}
\hfill
%-------------------------------------------------
% Panel (b): Resulting subforests poset
%-------------------------------------------------
\subfigure[]{
\begin{tikzpicture}[
    node distance=1.25cm and 2cm,
    every node/.style={
        draw,
        rounded corners,
        inner sep=2pt,
        font=\small
    },
    >=latex
]

\node (top) {$\rt[[[][]]]$};

\node (e1) [below left=of top] {$\rt[]\rt[[[]]]$};
\node (e2) [below=of top] {$\rt[]\rt[[[]]]$};
\node (e3) [below right=of top] {$\rt[]\rt[[][]]$};

\node (e12) [below left=of e2] {$\rt[]\rt[]\rt[[]]$};
\node (e13) [below right=of e2] {$\rt[]\rt[]\rt[[]]$};
\node (e23) [below=of e2] {$\rt[]\rt[]\rt[[]]$};

\node (all) [below=of e23] {$\rt[]\rt[]\rt[]\rt[]$};

\draw (top)--(e1);
\draw (top)--(e2);
\draw (top)--(e3);

\draw (e1)--(e12);
\draw (e1)--(e13);

\draw (e2)--(e12);
\draw (e2)--(e23);

\draw (e3)--(e13);
\draw (e3)--(e23);

\draw (e12)--(all);
\draw (e13)--(all);
\draw (e23)--(all);
\label{subfig:subforest_poset}
\end{tikzpicture}
}

\caption{A (colorless) rooted tree $t = \rt[[[][]]]$ and its associated edge-deletion poset. The elements of the posets correspond to subsets of edges removed from the original tree in \ref{subfig:subforest_deletion} and to subforests obtained after removing those edges in \ref{subfig:subforest_poset}. }
\label{fig:edge_poset}
\end{figure}

% Now define them on an arbitrary locally finite poset.
\subsection{Arbitrary locally finite posets}
The edge-deletion poset previously introduced  provides a concrete example of a finite partially ordered set. The algebraic structures associated with such posets are most naturally described in the more general setting of incidence algebras. We therefore recall the basic definitions and properties of incidence algebras over arbitrary locally finite posets before specializing them to the edge-deletion poset of a rooted tree.

Let $(P,\leq)$ be a locally finite partially ordered set, that is, a poset for which every interval
$$
[x,y]={z\in P:x\leq z\leq y}
$$
contains finitely many elements \cite{SCHMITT1994incidenceHopf}. The \textbf{incidence algebra} of $P$ over a field $\mathbb{K}$ is the vector space
$$
I(P)={f:P\times P\rightarrow\mathbb{K},\:\;f(x,y)=0\text{ whenever }x\not\leq y},
$$
equipped with pointwise addition and scalar multiplication. The multiplication is given by the convolution product
$$
(f*g)(x,y)
=
\sum_{x\leq z\leq y}
f(x,z),g(z,y),
$$
where the sum is finite by local finiteness of the poset \cite{aigner2012combinatorial,SCHMITT1994incidenceHopf}.

% Mobius inversion
\subsection{Distinguished elements of the incidence algebra}

Among the elements of the incidence algebra, two functions play a distinguished role. The first is the \emph{Kronecker delta},
\[
\delta(x,y)=
\begin{cases}
1,&x=y,\\
0,&x\neq y,
\end{cases}
\]
which is the identity element of the incidence algebra. Indeed, for every
$f\in I(P)$,
\[
f*\delta=\delta*f=f.
\]

The second is the \emph{zeta function},
\[
\zeta(x,y)=
\begin{cases}
1,&x\leq y,\\
0,&\text{otherwise},
\end{cases}
\]
which records the order relation of the poset \cite{aigner2012combinatorial}.

A fundamental result in the theory of incidence algebras states that the zeta
function is invertible with respect to the convolution product \cite{aigner2012combinatorial,rota1964foundations}. Its inverse,
denoted by $\mu$, is called the \emph{Möbius function} of the poset, that is,
\[
\mu*\zeta=\zeta*\mu=\delta.
\]

Equivalently, the Möbius function is uniquely determined by the recursive
relations $\mu(x,x)=1$, and, for $x<y$,
\[
\mu(x,y)
=
-\sum_{x\leq z<y}\mu(x,z),
\]
or, equivalently,
\[
\mu(x,y)
=
-\sum_{x<z\leq y}\mu(z,y).
\]

\subsection{Möbius inversion for Operator Splitting methods}

The Möbius function provides the inverse operation to summation over intervals.
This is formalized in the classical Möbius inversion theorem.

\begin{thm}[Möbius inversion, \cite{stanley2011enumerative}]\label{thm:mobius_inversion}
Let $f,g:P\rightarrow\mathbb{K}$ be functions on a locally finite poset
$(P,\leq)$. If
\[
g(x)=\sum_{y\leq x}f(y),
\]
for every $x\in P$, then
\[
f(x)=\sum_{y\leq x}\mu(y,x)\,g(y).
\]
Equivalently, if
\[
g(x)=\sum_{x\leq y}f(y),
\]
then
\[
f(x)=\sum_{x\leq y}\mu(x,y)\,g(y).
\]
\end{thm}
A particularly important example of Möbius inversion arises in the Boolean lattice $B_n(P)$, consisting of all subsets of an $n$-element set ordered by inclusion. In this case, the Möbius function is given explicitly by
\begin{equation}
\mu(A,B)=(-1)^{|B\setminus A|},\qquad A\subseteq B,
\end{equation}
so that the Möbius inversion reduces to the classical inclusion-exclusion principle \cite{SCHMITT1994incidenceHopf}. 

The edge-deletion poset associated with a rooted tree is precisely such a Boolean lattice: every subforest is uniquely determined by the subset of edges deleted from the original tree, and the partial order is induced by inclusion of these subsets. The following Proposition shows that, rather than viewing these identities as isolated combinatorial formulas, they are recognized as instances of Möbius inversion on the edge-deletion poset of a rooted tree.

Define, for the poset $P$, the lower sum operator $S_{\leq} $ on $ Map(P,K) $ given by $$ (S_{\leq}f)(x) := \sum_{y \leq x}f(y). $$ The fundamental theorem on Mobius inversion states that the difference operator $D_{\leq}$, inverse of $S_{\leq}$, exists and is completely determined by the Mobius function \cite{aigner2012combinatorial}. The lower difference operator is given by $$ (D_{\leq}g)(x) = \sum_{y \leq x} g(y) \mu(y,x). $$ Then, if $B(P)$ is the Boolean incidence algebra over an $n$-set $S$, then for any $f,g \in Map(B(P),K)$, with $K$ a field, we have \begin{equation} g(B) = \sum_{A \subseteq B}f(A) \iff f(B) = \sum_{A \subseteq B} (-1)^{|A\setminus B|}g(A). \end{equation} 

The immediate application to the theory of OS methods is the following proposition.

\begin{prop}
Let $t$ be a rooted tree, and let $\phi$ and $\phi^C$ denote the characters introduced in Section \ref{sec:background}. Then
\begin{equation} 
(-1)^{|E(T)|}\phi(T) = \sum_{F\subseteq T}\mu(F,T) \phi^C(F) \iff \phi^C(T) = \sum_{F\subseteq T}\mu(F,T) (-1)^{|E(T)|} \phi(F), 
\label{eq:mobius_phi_inverse}
\end{equation} 
where $\mu(F,T) = (-1)^{|E(T)| - |E(F)|} $.

Furthermore,
\[
D_{\leq}\phi^C
=
(-1)^{|E|}\phi,
\]
or equivalently,
\[
S_{\leq}\!\left((-1)^{|E|}\phi\right)
=
\phi^C.
\]
\end{prop}

\begin{proof}
    It follows from comparing \eqref{eq:phi_subgraphs} with Theorem \ref{thm:mobius_inversion}.
\end{proof}

\section{Hopf algebra of graphs}
\label{sec:hopf_algebra}
\subsection{Schmitt's incidence Hopf algebras}

Although incidence algebras are naturally associated with locally finite posets, a richer algebraic structure arises by considering hereditary families of finite intervals. This construction, introduced by Schmitt \cite{SCHMITT1994incidenceHopf,SCHMITT1987antipodes}, unifies numerous combinatorial Hopf algebras, including the Hopf algebra of graphs considered in this work.

Let $\mathcal P$ be a hereditary family of finite intervals, that is, a family of finite posets with unique minimal and maximal elements, closed under the formation of subintervals and Cartesian products. Suppose furthermore that the intervals in $\mathcal P$ are equipped with an order-compatible equivalence relation $\sim$, and denote by $[P]$ the equivalence class of an interval $P\in\mathcal P$. The vector space
\[
\mathcal H(\mathcal P)
=
\operatorname{span}_{\mathbb K}
\{[P]:P\in\mathcal P\}
\]
is endowed with the product
\[
[P_1]\cdot[P_2]
=
[P_1\times P_2],
\]
where $P_1\times P_2$ denotes the Cartesian product of posets.

The coalgebra structure is defined by decomposing an interval at every intermediate element \cite{SCHMITT1987antipodes,SCHMITT1994incidenceHopf}. More precisely,
\[
\Delta([P])
=
\sum_{x\in P}
[\hat 0,x]\otimes[x,\hat 1],
\]
where $\hat 0$ and $\hat 1$ denote the minimal and maximal elements of $P$. The counit is given by
\[
\varepsilon([P])
=
\begin{cases}
1,& |P|=1,\\
0,&\text{otherwise}.
\end{cases}
\]

The coproduct is coassociative, while the product is associative. Under the compatibility assumptions introduced by Schmitt, these structures satisfy the bialgebra axioms \cite{schmitt1993hopf}. Moreover, if $\mathcal H(\mathcal P)$ is connected and graded, then it admits a unique antipode, thereby becoming a Hopf algebra \cite{milnor1965structure}. The antipode is given recursively in \cite{SCHMITT1994incidenceHopf} by
\[
S([P])
=
-
[P]
-
\sum_{\hat 0<x<\hat 1}
S([\hat 0,x])\,[x,\hat 1].
\]

\subsection{The Hopf algebra of graphs}

The Hopf algebra of graphs was introduced by Schmitt as an instance of an incidence Hopf algebra associated with a hereditary family of graphs \cite{SCHMITT1994incidenceHopf}. We briefly recall its construction, paying particular attention to the antipode.

Let $\mathcal G$ denote the family of finite simple graphs, considered up to isomorphism. The associated Hopf algebra is
\[
\mathcal H(\mathcal G)
=
\operatorname{span}_{\mathbb K}
\{[G]:G\in\mathcal G\},
\]
whose basis consists of the isomorphism classes of finite graphs.

The product is induced by disjoint union,
\[
[G_1]\cdot[G_2]
=
[G_1\sqcup G_2],
\]
while the coproduct is given by
\[
\Delta([G])
=
\sum_{U\subseteq V(G)}
[G|_U]
\otimes
[G|_{V(G)\setminus U}],
\]
where $G|_U$ denotes the induced subgraph on the vertex set $U$.

A formula for the antipode in function of the sets of acyclic orientations was introduced in \cite{humpert2012incidence}. Defining an acyclic orientation of a graph $G$ as a choice of orientation of all the edges that admits no directed cycles, they introduce $\mathcal{A}(G)$ as the set of acyclic orientations of $G$, $\mathcal{F}(G)$ as the set of flats of $G$, and define $a(G) := |\mathcal{A}(G)|$. Then, the combinatorial antipode formula of Humpert et al. for a graph $G= (V,E)$, based on Takeuchi’s formula
for connected Hopf algebras, reads 
\begin{equation}
    S(G) = \sum_{F\in \mathcal{F}(G)}(-1)^{|V| - rk(F)}a(G/F)G_{V,F},
    \label{eq:antipode_graphs}
\end{equation}
where the rank $rk(F)$, for $F\subseteq E(G)$, is the size of any maximal acyclic
subset of $F$.

\subsection{The antipode and the group of characters}\label{sec:antipode}
\begin{thm}\label{thm:K_antipode}
    For a rooted tree $t$,
\begin{equation}
S(t)=(-1)^{|V(t)|}K(t),
\end{equation}
where $S$ is the antipode of Schmitt's Hopf algebra of graphs.
\end{thm}
\begin{proof}
    For a rooted tree forest, $a(t) = 1$. Consequently, up to the parity factor $(-1)^{|V(T)|}$, the operator $K$ in \eqref{eq:K_expansion} coincides with the antipode \eqref{eq:antipode_graphs} of the Hopf algebra of graphs.
\end{proof}

We focus on the restriction of $\mathcal{G}$ to the subalgebra of graphs corresponding to rooted trees, $\mathcal H(\mathcal T)$, graded by the number of vertices.
\begin{corollary}
Under the identification of rooted trees with the Boolean lattice of their
subforests, the antipode of $\mathcal{H(\mathcal{G})}$ satisfies
\[
S=(-1)^{|V|}D_{\leq},
\]
where $D_{\leq}$ denotes the Möbius difference operator.
\end{corollary}

Note that the action of the operator $K$ and the properties of the characters $\phi,\phi^C$ can be studied for rooted trees with any coloring of their vertices, since the terms obtained after applying $K$ to a rooted tree depends only on the set of edges, and the tree factorial is calculated only on the shape of the tree and not on the labels of its vertices.

Before digging into the applications of Theorem \ref{thm:K_antipode} to OS methods and studying the characters on $K$, we show that Proposition \ref{prop:involution} can be derived, directly and without heavy computation, as a Corollary of Theorem \ref{thm:K_antipode} by noting that the antipode of an incidence Hopf algebra is an involution \cite{SCHMITT1994incidenceHopf}.
\begin{corollary}
    $K$ is an involution:
    \begin{equation}
        K^2(t) := (K\circ K)(t)  = t.
    \end{equation}
\end{corollary}

Now, we focus on the Characters Group $\mathcal{G}(\mathcal H(\mathcal T),\mathbb{K})$. The elements of $\mathcal{G}(\mathcal H(\mathcal T),\mathbb{K})$ form a group with product given by convolution \cite{bogfjellmo2019characters}
\begin{equation}
    \phi*\varphi = m\circ (\phi\otimes \varphi)\circ \Delta, \quad \phi,\varphi \in \mathcal{G}(\mathcal{H},\mathbb{K}),
\end{equation}
where $m:\mathcal H(\mathcal T)\times \mathcal H(\mathcal T) \rightarrow \mathcal H(\mathcal T)$ is the multiplication and $\Delta$ is the coproduct. For example, the character group of the Hopf algebra of rooted trees of Connes and Kreimer is known as the \textit{Butcher group} in numerical analysis \cite{butcher1972algebraic2}.

For all $a,b \in \mathcal{H}$, the characters $\varphi: \mathcal{H}\rightarrow \mathbb{H}$ satisfy multiplicativity: $\varphi(ab) = \varphi(a)\varphi(b)$, and $\varphi(1) = 1$. In addition, the inverse in convolution of a given character $\varphi$ is obtained via the antipode of the Hopf algebra: $\varphi^{-1}=\varphi \circ S$. For a character $\varphi$, there is a canonical involution characters depending on the grading $n$ of the Hopf algebra defined by $\bar{\varphi} := (-1)^n\varphi$ \cite{Aguiar_Bergeron_Sottile_2006}.

\begin{thm}\label{thm:phi_inverse}
Let $\phi,\phi^C\in\mathcal G(\mathcal H(\mathcal T),\mathbb K)$
denote the characters associated with the elementary weights
defined recursively in \eqref{eq:phi_def} and
\eqref{eq:phi_complementary}. Then, $\phi,\phi^C \in \mathcal{G}(\mathcal H(\mathcal T),\mathbb{K})$, and
    \begin{align}
        \phi^{-1} &= \bar {\phi}^C,\\
        (\phi^C)^{-1} &= \bar {\phi}.
    \end{align}
\end{thm}
\begin{proof}
    Since we assume first-order for the coefficients of the OS methods, the characters satisfy $\phi(t) = \phi^C(t) = 1$ if $|t| =1$. Then,
$$
\phi^{-1} = \phi\circ S = (-1)^{|V|}\phi \circ K= (-1)^{|V|}\phi^C.
$$
\end{proof}
Theorem \ref{thm:phi_inverse} provides a natural Hopf algebraic interpretation of the
characters $\phi$ and $\phi^C$. Rather than being two independent characters,
they are related through the composition of two canonical operations on the
character group: convolution inversion and the canonical involution. Consequently, the characters $\phi$ and
$\phi^C$ form a two-element orbit under this natural action on the group of
characters.

Now, with the theory of posets, Theorem \ref{thm:phi_homeormorphic_candycanes} can be understood from the perspective of combinatorial Hopf algebras: since two rooted trees belonging to the same equivalence class as unrooted trees possess the same underlying graph and therefore the same collection of edges, the edge-deletion poset depends only on the subsets of edges removed. Then, every subforest in the edge-deletion poset of one tree is naturally associated with the subforest obtained by retaining the same subset of edges and rerooting each connected component according to the new choice of root. Although the rooted structure of the individual components may change, the inclusion relations between edge sets are preserved. This correspondence defines an order isomorphism between the two edge-deletion posets, explaining why the combinatorial identities arising from the antipode, and hence the tree-factorial identities, are the same for both rooted trees.

%\gabriel{Should I keep this or not?}
%On the other hand, the relation between combinatorial Hopf algebras and chromatic polynomials has been noted, for example, in \cite{SCHMITT1994incidenceHopf,Aguiar_Bergeron_Sottile_2006,humpert2012incidence}. The Euler characteristic of a poset can be obtained from the Mobius function:
% \begin{equation}
%     \chi(P) = \sum_{i=0}^n(-1)^if_i(P),
% \end{equation}
% where $n = rk(P)$ and $f_i(P)$ counts the number of elements of P of rank $i$. Since the operator $K$ coincides with the Möbius difference operator, by a factor of $(-1)^{|V|}$, on the Boolean lattice of subforests, composing it with the constant character recovers the counit.

%\begin{prop}
%If $\zeta \in \mathcal{G(\mathcal{H(G)},\mathbb{K})}$ is the character such that $\zeta(t) = 1 \; \forall t$, then
%    \begin{equation}
%    (\zeta \circ K)(t) = (-1)^{|V|}\varepsilon(t) =
%    \begin{cases}
%        1, \quad |t| = 1,\\
%        0, \quad else.
%    \end{cases}
%\end{equation}
%\end{prop}

\section{Applications to numerical analysis}
\label{sec:applications}
\subsection{Operator Splitting methods}
The notion of an adjoint method plays a fundamental role in geometric numerical
integration, particularly in the construction of symmetric and higher-order
composition methods \cite{hairer2003geometric}. Let $\Phi_h$ denote the numerical flow associated with a
one-step method with time step $h$. The \emph{adjoint method} is defined by
\begin{equation}
    \Phi_h^*:=\Phi_{-h}^{-1},
\end{equation}
that is, the method obtained by reversing the numerical flow. A method is said
to be \emph{symmetric} or \emph{self-adjoint} if $\Phi_h=\Phi_h^*$.

For OS methods, the adjoint admits a particularly simple
description. If
\[
\Phi_h
=
\varphi_{\alpha_1^\mu h}^{\mu}
\varphi_{\alpha_1^\nu h}^{\nu}
\cdots
\varphi_{\alpha_s^\mu h}^{\mu}
\varphi_{\alpha_s^\nu h}^{\nu},
\]
where $\varphi_h^{(i)}$ denotes the flow associated with the $i$-th operator,
then
\[
\Phi_h^*
=
\varphi_{\alpha_s^\nu h}^{\nu}
\varphi_{\alpha_s^\mu h}^{\mu}
\cdots
\varphi_{\alpha_1^\nu h}^{\nu}
\varphi_{\alpha_1^\mu h}^{\mu},
\]
that is, the adjoint is obtained by reversing the order of the composition.
Consequently, the elementary weights of the adjoint method are obtained by
reversing the OS coefficients together with the sign change
induced by the substitution $h\mapsto -h$. 

The following theorem is probably the main result of this section as it illustrates the potential of combinatorial Hopf algebras for geometric integration, in particular for OS methods.

\begin{corollary}\label{thm:adjoint_character}
Let $\Phi$ be an FSRK method with character $\phi$. Let $\Phi^*$ denote its
adjoint method, and let $\phi^*$ be the corresponding character. Then
\[
\phi^{-1}=\overline{\phi^C} = \phi^*.
\]
Equivalently, the inverse in convolution of the character coincides with the elementary weights of the adjoint method.
\end{corollary}
\begin{proof}
By definition, the elementary weights of the adjoint method are obtained by reversing the order of the OS coefficients,
\[
\alpha_i^\ell \longmapsto \alpha_{s+1-i}^\ell,
\qquad
\ell\in\{\mu,\nu\},
\]
together with the sign change induced by the substitution $h\mapsto -h$. The index reversal
\(
i\mapsto s+1-i
\)
is precisely the transformation defining the complementary character $\phi^C$, while the sign change is accounted for by the canonical involution. Hence,
\[
\phi^*=\overline{\phi^C}.
\]
\end{proof}
\noindent
Now, a method is symmetric if and only if \cite{hairer2003geometric}
\[
\alpha_i^\ell=\alpha_{s+1-i}^\ell,
\qquad
i=1,\ldots,s,\quad
\ell\in\{\mu,\nu\},
\]
that is, if the OS coefficients form a palindrome. 
\begin{corollary}
    If an OS method is symmetric, then $\phi = \phi^C$.
\end{corollary}
The following result is already found in the literature (see \cite{hairer2003geometric}, Chapter II), and we rederive it using the theory of combinatorial Hopf algebras.
\begin{corollary}[Symmetry Theorem]
Let $\Phi_h$ be a symmetric OS method. Then its local error
expansion contains only odd powers of the step size $h$. Consequently, the
order of the method is necessarily even.
\end{corollary}

\subsection{Partitioned Runge-Kutta methods}
\label{sec:partitioned}
The combinatorial ideas developed in this work are not restricted to OS methods. As an illustration, consider the identity of Corollary \ref{corollary:edge_rooted}, which appears in the order theory of partitioned Runge--Kutta methods \cite{Abia1993partitioned}. 

Partitioned Runge--Kutta (PRK) methods are a generalization of classical RK methods designed for systems whose variables naturally split into different components. Instead of advancing every component with the same Butcher tableau, PRK methods employ two coupled RK schemes, allowing each partition to be integrated according to its own coefficients \cite{hairer2003geometric}.

Here we use the Hopf algebra of graphs to provide a combinatorial explanation of Corollary \ref{corollary:edge_rooted}. The connection between tree factorials and partially ordered sets is provided by the following classical result of Knuth.

\begin{thm}[Knuth,\cite{Knuth1998programming}]
Let $T_r$ be a rooted tree with $n$ vertices. Then the number of linear extensions
of the associated tree poset is
\[
|\mathcal{L}(T_r)|=\frac{n!}{T_r!},
\]
where $\mathcal{L}(T_r)$ denotes the set of linear extensions of $T_r$ .
\end{thm}

\begin{prop}\label{prop:linear_extension}
Let $T$ be a tree and let $e=(r,r')$ be an edge of $T$. Let
$F=c_e(T)$ be the rooted forest obtained by deleting $e$ and rooting the
resulting connected components at $r$ and $r'$, respectively. Then the set of
linear extensions of $F$ admits the disjoint decomposition
\[
\mathcal{L}(F)
=
\mathcal{L}_{r<r'}(F)
\sqcup
\mathcal{L}_{r'<r}(F),
\]
where
\[
\mathcal{L}_{r<r'}(F)
=
\{\sigma\in\mathcal{L}(F):\sigma(r)<\sigma(r')\},
\]
and
\[
\mathcal{L}_{r'<r}(F)
=
\{\sigma\in\mathcal{L}(F):\sigma(r')<\sigma(r)\}.
\]
Moreover, there are natural bijections
\[
\mathcal{L}_{r<r'}(F)\cong\mathcal{L}(T_r),
\qquad
\mathcal{L}_{r'<r}(F)\cong\mathcal{L}(T_{r'}).
\]
\end{prop}
\begin{proof}
Every linear extension of the forest poset $F$ satisfies exactly one of the
mutually exclusive relations $r<r'\qquad$ or $\qquad r'<r$.
Hence,
\[
\mathcal{L}(F)
=
\mathcal{L}_{r<r'}(F)
\sqcup
\mathcal{L}_{r'<r}(F).
\]

Suppose that $\sigma\in\mathcal{L}_{r<r'}(F)$. Adding the relation $r<r'$
to the poset $F$ reconstructs the rooted tree poset $T_r$. Since $\sigma$
already satisfies $r<r'$, it follows that $\sigma$ is a linear extension of
$T_r$. Conversely, every linear extension of $T_r$ is also a linear extension
of $F$ satisfying $r<r'$, since removing the relation $r<r'$ cannot violate
any of the remaining order relations. Therefore,
\[
\mathcal{L}_{r<r'}(F)\cong\mathcal{L}(T_r).
\]

The same argument with the roles of $r$ and $r'$ interchanged yields the
bijection $\mathcal{L}_{r'<r}(F)\cong\mathcal{L}(T_{r'})$.
\end{proof}

\noindent
Taking cardinalities and using the bijections of Proposition~\ref{prop:linear_extension} gives
\[
|\mathcal{L}(F)|
=
|\mathcal{L}_{r<r'}(F)|
+
|\mathcal{L}_{r'<r}(F)|
=
|\mathcal{L}(T_r)|
+
|\mathcal{L}(T_{r'})|.
\]
Since all three posets have the same underlying vertex set, Knuth's theorem
yields
\[
|\mathcal{L}(F)|
=
\frac{|T|!}{F!},
\qquad
|\mathcal{L}(T_r)|
=
\frac{|T|!}{T_r!},
\qquad
|\mathcal{L}(T_{r'})|
=
\frac{|T|!}{T_{r'}!},
\]
where we have used the convention $(F=t_1\cdot t_2)! = t_1!\,t_2!$. Therefore,
\[
\frac{1}{F!}
=
\frac{1}{T_r!}
+
\frac{1}{T_{r'}!},
\]
which is precisely the rerooting identity of Corollary~\ref{corollary:edge_rooted}. This decomposition of linear extensions is therefore a manifestation of the combinatorial structure encoded by the edge-deletion poset, providing a natural explanation of the rerooting identity.

\section{Conclusions}
Motivated by the identity that expresses the elementary weight of a rooted tree as an inclusion-exclusion expansion over all its rooted subforests, we established a connection between this expansion and the antipode of a combinatorial Hopf algebra. Theorem~\ref{thm:K_antipode} provides the principal algebraic interpretation of the present work. Although the operator $K$ was originally introduced in the study of FSRK methods to describe the dependence of elementary weights on subforests, the theorem shows that it is, in fact, a fundamental object in the Hopf algebraic structure of graphs. Consequently, many of the combinatorial identities satisfied by $K$ are no longer isolated phenomena, but rather manifestations of general properties of the antipode. This identification establishes a direct connection between the combinatorics of OS methods and the theory of incidence Hopf algebras.

The edge-deletion poset provides a complementary perspective. Its incidence algebra and Möbius function encode the same alternating-sign structure, linking the subforest expansions to Möbius inversion.
Our motivation stems from the study of bipartite
rooted trees, whose coloring properties naturally relate them to the chromatic
theory of graphs \cite{Aguiar_Bergeron_Sottile_2006}. Indeed, a graph is bipartite if and only if it admits a proper
$2$-coloring, while the chromatic polynomial is classically expressed as a
Möbius inversion over the Boolean lattice of edge subsets. In the present context, it provides a systematic way for interpreting the identities associated with the operator $K$.

Beyond providing a unified framework for previously known identities in OS and PRK methods, such as Theorems~\ref{thm:phiC} and~\ref{thm:phi_homeormorphic_candycanes} and the identity~\eqref{eq:partitioned}, this approach has enabled us to derive new combinatorial identities, including Propositions~\ref{prop:involution} and~\ref{prop:degree}. 
In particular, adjoint operator splitting methods admit an algebraic interpretation involving Möbius inversion and the canonical involution of the original method.  Furthermore, just as Runge--Kutta methods can be understood through the structure of the Connes--Kreimer--Butcher Hopf algebra, OS methods can be naturally described within the framework of Schmitt's Hopf algebra of graphs.
This connection opens a route to studying further algebraic and combinatorial properties of geometric numerical integrators.

\bibliographystyle{plain}
\bibliography{bib}
\end{document}